\documentclass[10pt,onecolumn]{IEEEtran}       
\usepackage{amsfonts,amssymb,amsmath}
\usepackage{indentfirst,latexsym,amscd,algorithmic}
\usepackage{amsbsy,amsthm,caption}
\usepackage[para]{threeparttable}
\usepackage{graphicx,multirow,bm}
\usepackage{caption}
\usepackage{color}
\usepackage{array}
\usepackage{subfigure}
\usepackage[ruled,vlined,commentsnumbered]{algorithm2e}
\usepackage{pifont,marvosym}
\usepackage{tikz}
\usepackage{longtable}
\usepackage{pdfpages}
\usepackage{geometry} 
\usepackage{setspace}
\usepackage[justification=centering]{caption}
\usepackage[noadjust]{cite}
\renewcommand{\vec}[1]{\boldsymbol{#1}} 

\newcommand{\C}{{\cal C}}
\newcommand{\F}{{\cal F}}
\newcommand{\A}{{\cal A}}

\theoremstyle{mythm}
\newtheorem{definition}{Definition}[section]
\newtheorem{lemma}{Lemma}[section]
\newtheorem{theorem}{Theorem}[section]
\newtheorem{corollary}[theorem]{Corollary}
\newtheorem{remark}{Remark}
\newtheorem{example}{Example}

\begin{document}

\title{Strongly separable matrices for nonadaptive combinatorial group testing
\author{Jinping Fan, Hung-Lin Fu, Yujie Gu, Ying Miao, and Maiko Shigeno}
\thanks{J. Fan is with the Department of Policy and Planning Sciences, Graduate School of Systems and Information Engineering, University of Tsukuba, Tsukuba, Ibaraki 305-8573, Japan (e-mail: j.fan.math@gmail.com).}
\thanks{H.-L. Fu is with the Department of Applied Mathematics, National Chiao Tung University, Hsinchu 30010, Taiwan (e-mail: hlfu@math.nctu.edu.tw).}
\thanks{Y. Gu is with the Department of Electrical Engineering-Systems, Tel Aviv University, Tel Aviv 6997801, Israel (e-mail: guyujie2016@gmail.com).}
\thanks{Y. Miao and M. Shigeno are with the Faculty of Engineering, Information and Systems, University of Tsukuba, Tsukuba, Ibaraki 305-8573, Japan (e-mails: miao@sk.tsukuba.ac.jp; maiko@sk.tsukuba.ac.jp).}
}

\date{}

\maketitle

\begin{abstract}
In nonadaptive combinatorial group testing (CGT), it is desirable to identify a small set of up to $d$ defectives from a large population of $n$ items with as few tests (i.e. large rate) and efficient identifying algorithm as possible. In the literature, $d$-disjunct matrices ($d$-DM) and $\bar{d}$-separable matrices ($\bar{d}$-SM) are two classical combinatorial structures having been studied for several decades. It is well-known that a $d$-DM provides a more efficient identifying algorithm than a $\bar{d}$-SM, while a $\bar{d}$-SM could have a larger rate than a $d$-DM. In order to combine the advantages of these two structures, in this paper, we introduce a new notion of \emph{strongly $d$-separable matrix} ($d$-SSM) for nonadaptive CGT and show that a $d$-SSM has the same identifying ability as a $d$-DM, but much weaker requirements than a $d$-DM. Accordingly, the general bounds on the largest rate of a $d$-SSM are established. Moreover, by the random coding method with expurgation, we derive an improved lower bound on the largest rate of a $2$-SSM which is much higher than the best known result of a $2$-DM.
\end{abstract}


\section{Introduction}
\label{sec:1}

Group testing was introduced by Dorfman \cite{RD} in 1940s for a large-scaled blood testing program. The object of this program was testing a large number of blood samples to determine the \textit{defective} (or \textit{positive}) ones. Instead of testing one by one, group testing was proposed to pool all the blood samples into groups and perform a test to each group. If the testing outcome of a group is \textit{positive}, it means that at least one defective is contained in this group. If the testing outcome of a group is \textit{negative}, then this group is considered containing no defective samples. In general, there are two types of group testing models. In an \textit{adaptive} (or \textit{sequential}) setting, the group arrangement of the next test is always based on the previous testing outcomes. In a \textit{nonadaptive} setting, all the group arrangements are designed in advance and all the tests are conducted simultaneously. On the other hand, group testing can be roughly divided into two categories: \textit{probabilistic group testing} (PGT) and \textit{combinatorial group testing} (CGT). In PGT, the defective samples are assumed to follow some probability distribution, while in CGT, the number of defective samples is usually assumed to be no more than a fixed positive integer. In this paper, we focus on \textit{nonadaptive CGT} (or \textit{pooling design}) in the noiseless scenario, which has been studied extensively due to its applications in a variety of fields such as DNA library screening, network security, multi-access communication and so on, see \cite{DH,DH2006,D,KS,ND,W} for example.

A nonadaptive CGT scheme can be represented by a binary matrix whose rows are indexed by the groups and columns are indexed by the items to be tested. Suppose that there are $n$ items among which at most $d$ $(\ll n)$ are positive. Let $M$ be the $t\times n$ testing matrix where the entry $M_j(i)=1$ if the $j$th item is contained in the $i$th test and $M_j(i)=0$ otherwise. Note that each column of $M$ corresponds to an item and each row corresponds to a test. The result of a test is $1$ (positive) if the test contains at least one positive item and $0$ (negative) otherwise. After performing all the $t$ tests simultaneously, one could observe the testing outcome $\vec{r}=(\vec{r}(1),\ldots,\vec{r}(t))\in\{0,1\}^t$ where $\vec{r}(i)$ is the result of the $i$th test. It is easily seen that $\vec{r}$ is the \textit{Boolean sum} of the column vectors of $M$ indexed by all the positive items. The problem of studying nonadaptive CGT is to design the testing matrix $M$ such that all the positives could be identified based on $M$ and the testing outcome $\vec{r}$. The goal is to decrease the computational complexity of the identifying algorithm and to minimize the number of tests needed given the number of items to be tested, or equivalently, to maximize the number of items to be tested given the number of tests, or in other words, to explore the \textit{largest rate} of the testing matrix $M$.

In the literature, \textit{disjunct matrices} and \textit{separable matrices} are two classical combinatorial structures for nonadaptive CGT. Disjunct matrices and separable matrices were first studied by Kautz and Singleton \cite{KS} for file retrieval under the name of \textit{superimposed codes}, and later were extensively investigated under the concepts of \textit{cover-free families} and \textit{union-free families} respectively, see \cite{EFF,EFF1985,FF} for example. The definitions of disjunct matrices and separable matrices could be found in \cite{DH}.
\begin{definition}{\rm (\cite{DH})}\label{dm-sm}
Let $n,t,d\geq 2$ be integers and $M$ be a binary matrix of size $t\times n$.
\begin{enumerate}
  \item $M$ is called a \textit{$d$-disjunct matrix}, or briefly $d$-DM, if the Boolean sum of any $d$ column vectors of $M$ does not cover any other one.
  \item $M$ is called a \textit{$\bar{d}$-separable matrix}, or briefly $\bar{d}$-SM, if the Boolean sums of $\le d$ column vectors of $M$ are all distinct.
\end{enumerate}
\end{definition}

It was shown in \cite{DH} that a $d$-DM, as well as a $\bar{d}$-SM, could be utilized in nonadaptive CGT to identify any set of positives with size no more than $d$, but both have their advantages and disadvantages. In general, the computational complexity of the identifying algorithm based on a $d$-DM of size $t\times n$ is $O(tn)$, while that based on a $\bar{d}$-SM of size $t\times n$ is $O(tn^d)$, explicating that a $d$-DM provides a more efficient identifying algorithm than a $\bar{d}$-SM. However, given the number of rows, a $\bar{d}$-SM contains more columns, implying a higher rate, than a $d$-DM. In order to combine the advantages of these two structures, in this paper, we introduce a new notion of \textit{strongly $d$-separable matrix} ($d$-SSM) for nonadaptive CGT which has weaker requirements than $d$-DM but the same identifying ability as $d$-DM. It is also shown that a $d$-SSM has a larger rate than a $d$-DM.

The rest of this paper is organized as follows. In Section \ref{sec:2}, we give the definition of SSM and establish the relationships among SSM, DM and SM. We show that a $d$-SSM could identify any set of up to $d$ positives as efficiently as a $d$-DM. In Section \ref{sec:4}, we first give the general bounds on the largest rate of a $d$-SSM from the known results on DM and SM, and then derive an improved lower bound on the largest rate of a $2$-SSM by the random coding method with expurgation. The conclusion is drawn in Section \ref{sec:5}.

\section{Strongly separable matrices}
\label{sec:2}

In this section, we first introduce the notion of $d$-SSM and investigate the relationships among DM, SSM and SM. Next we provide an identifying algorithm based on a $d$-SSM and prove that a $d$-SSM has the same identifying ability as a $d$-DM.

Let $n,t,d$ be integers with $n\ge d\ge 2$, and $M$ be a binary matrix of size $t\times n$. Denote $[t]=\{1,\ldots,t\}$ and $[n]=\{1,\ldots,n\}$. Let $\F=\{\vec{c}_1,\dots,\vec{c}_n\}\subseteq\{0,1\}^t$ be the set of column vectors of $M$ where $\vec{c}_j=(\vec{c}_j(1),\ldots,\vec{c}_j(t))\in\{0,1\}^t$ for any $j\in[n]$. We say a vector $\vec{c}_j$ \textit{covers} a vector $\vec{c}_k$ if for any $i\in [t]$, $\vec{c}_k(i)=1$ implies $\vec{c}_j(i)=1$.

\begin{definition}\label{ssm}
A $t\times n$ binary matrix $M$ is called a \textit{strongly $d$-separable matrix}, or briefly $d$-{\rm SSM}, if for any $\F_0\subseteq\F$ with $|\F_0|=d$, we have
\begin{equation}\label{ssm-1}
  \bigcap\limits_{\F'\in \mathcal{U}(\F_0)}\F'=\F_0,
\end{equation}
where
\begin{equation}\label{ssm-2}
  \mathcal{U}(\F_0)=\Big\{\F'\subseteq\F: \bigvee_{\vec{c}\in\F_0}\vec{c}=\bigvee_{\vec{c}\in\F'}\vec{c}\Big\}.
\end{equation}
\end{definition}

\begin{remark}\label{ssm:2}
An equivalent description of Definition \ref{ssm} is as follows. A $t\times n$ binary matrix $M$ is a $d$-SSM if for any $\F_0,\F'\subseteq\F$ with $|\F_0|=d$,
$\F'\in\mathcal{U}(\F_0)$ implies that $\F_0\subseteq\F'$.
\end{remark}

\begin{remark}\label{ssm:1}
We call $M$ a $\bar{d}$-SSM if the condition $|\F_0|=d$ in Definition \ref{ssm} is replaced by $1\le|\F_0|\le d$.
\end{remark}

It is obvious that a $\bar{d}$-SSM is always a $d$-SSM. The following observation shows that the converse also holds.

\begin{lemma}\label{equivalence}
A $t\times n$ matrix $M$ is a $d$-SSM if and only if $M$ is a $\bar{d}$-SSM.
\end{lemma}
\begin{IEEEproof}
It is enough to show the necessity. Assume that $M$ is a $d$-SSM of size $t\times n$, but not a $\bar{d}$-SSM. Let $\F=\{\vec{c}_1,\dots,\vec{c}_n\}$ be the set of column vectors of $M$. By Definition \ref{ssm} and Remarks \ref{ssm:2} and \ref{ssm:1}, there must exist $\F_0,\F'\subseteq\F$ with $1\le|\F_0|<d$ such that $\F'\in \mathcal{U}(\F_0)$ but $\F_0\not\subseteq\F'$. Suppose that $|\F_0|=j$, $1\le j\le d-1$. If $\F'\setminus\F_0=\emptyset$, there exists $\F_1\subseteq\F\setminus\F_0$ such that $|\F_0\cup\F_1|=d$ since $n\ge d$. Then we have $\F'\cup\F_1\in \mathcal{U}(\F_0\cup\F_1)$ but $\F_0\cup\F_1\not\subseteq \F'\cup\F_1$, a contradiction to the definition of $d$-SSM. Thus, $\F'\setminus\F_0\neq\emptyset$ and $\vec{c}'$ is covered by $\bigvee_{\vec{c}\in\F_0}\vec{c}$ for any $\vec{c}'\in\F'\setminus\F_0$. We discuss $|\F'\setminus\F_0|$ based on the following two cases.
\begin{enumerate}[1)]
  \item If $|\F'\setminus\F_0|\ge d-j$, then there exists $\F_1\subseteq\F'\setminus\F_0$ such that $|\F_0\cup\F_1|=d$ and $\F'\in \mathcal{U}(\F_0\cup\F_1)$. Since $\F_0\not\subseteq\F'$, we have $\F_0\cup\F_1\not\subseteq \F'$, a contradiction to the definition of $d$-SSM.

  \item If $|\F'\setminus\F_0|< d-j$, then $|\F_0\cup\F'|<d$. Since $n\ge d$, there exists $\F_2\subseteq\F\setminus(\F_0\cup\F')$ such that $|\F_0\cup\F'\cup\F_2|=d$ and $\F'\cup\F_2\in \mathcal{U}(\F_0\cup\F'\cup\F_2)$, but $\F_0\cup\F'\cup\F_2\not\subseteq\F'\cup\F_2$, also a contradiction.
\end{enumerate}

The conclusion follows.
\end{IEEEproof}

The relationship between DM and SM was investigated in \cite{CH,DH}.

\begin{lemma}{\rm{(\cite{CH,DH})}}\label{relationship of DM and SM}
A $d$-DM is a $\bar{d}$-SM and a $\bar{d}$-SM is a $(d-1)$-DM.
\end{lemma}

The following lemma shows that a $d$-SSM lies between a $d$-DM and a $\bar{d}$-SM.

\begin{lemma}\label{relationship2}
A $d$-DM is a $d$-SSM and a $d$-SSM is a $\bar{d}$-SM.
\end{lemma}
\begin{IEEEproof}
We first show that a $d$-DM is a $\bar{d}$-SSM. Let $A$ be a $d$-DM and $\A$ be the set the column vectors of $A$. By Definition \ref{dm-sm}, for any $\A_0\subseteq\A$ with $|\A_0|=d$, we have $\mathcal{U}(\A_0)=\{\A_0\}$ where $\mathcal{U}(\A_0)$ is defined by (\ref{ssm-2}). Then we obtain $\bigcap_{\A'\in \mathcal{U}(\A_0)}\A'=\A_0$ implying that $A$ is a $d$-SSM according to Definition \ref{ssm}.

Now we prove that a $d$-SSM is a $\bar{d}$-SM. Let $M$ be a $d$-SSM and $\F$ be the set the column vectors of $M$. If $M$ is not a $\bar{d}$-SM, then by Definition \ref{dm-sm}, there exist distinct $\F_1,\F_2\subseteq\F$ with $1\le |\F_1|,|\F_2|\le d$ such that $\bigvee_{\vec{c}\in\F_1}\vec{c}=\bigvee_{\vec{c}\in\F_2}\vec{c}$. By (\ref{ssm-2}), we have $\F_1\in \mathcal{U}(\F_2)$ and $\F_2\in \mathcal{U}(\F_1)$. Since $M$ is a $d$-SSM, by Lemma \ref{equivalence} we have $\F_1\subseteq\F_2$ and $\F_2\subseteq\F_1$ which implies $\F_1=\F_2$, a contradiction to the assumption.
\end{IEEEproof}

From Lemma \ref{relationship2} we know that a $d$-DM is always a $d$-SSM. But a $d$-SSM might not be a $d$-DM. To show this, we give an example below.

\begin{example}
Let $M$ be a binary matrix of size $7\times 8$ as defined below. 
\begin{equation*}
\begin{aligned}
M & =\begin{matrix}
\vec{c}_1 \hspace{0.16cm} \vec{c}_2 \hspace{0.175cm} \vec{c}_3 \hspace{0.17cm} \vec{c}_4 \hspace{0.18cm} \vec{c}_5 \hspace{0.2cm} \vec{c}_6 \hspace{0.18cm} \vec{c}_7 \hspace{0.18cm} \vec{c}_8 \\
\begin{bmatrix}
     1 & 0 & 0 & 0 & 0 & 0 & 0 & 1 \\
     1 & 1 & 0 & 0 & 0 & 0 & 0 & 0 \\
     0 & 1 & 1 & 0 & 0 & 1 & 0 & 0 \\
     0 & 0 & 1 & 1 & 0 & 0 & 0 & 0 \\
     0 & 0 & 0 & 1 & 0 & 0 & 1 & 0 \\
     0 & 0 & 0 & 0 & 1 & 1 & 0 & 0 \\
     0 & 0 & 0 & 0 & 1 & 0 & 1 & 1
\end{bmatrix}
\end{matrix}
\end{aligned}
\end{equation*}
It is easy to check that $M$ is a $2$-SSM and thus a $\bar{2}$-SM by Lemma \ref{relationship2}. But $M$ is not a $2$-DM since $\vec{c}_2$ is covered by $\vec{c}_1\vee\vec{c}_3$.
\end{example}

By Lemma \ref{relationship2}, a $d$-SSM, as well as a $\bar{d}$-SM, has weaker requirements than a $d$-DM. However, we will prove that a $d$-SSM could determine any set of positives with size no more than $d$ as efficiently as a $d$-DM.

\begin{theorem}\label{complexity}
A $t\times n$ $d$-SSM could identify any set of $\le d$ positives among $n$ items with $t$ tests by applying Algorithm \ref{algorithm 1}, and the computational complexity of Algorithm \ref{algorithm 1} is $O(tn)$.
\end{theorem}
\begin{IEEEproof}
Denote $[n]$ as $n$ items to be tested. Let $M$ be the testing matrix which is a $d$-SSM of size $t\times n$, and $\F=\{\vec{c}_1,\vec{c}_2,\ldots,\vec{c}_n\}$ be the set the column vectors of $M$. Suppose that $P_0\subseteq[n]$ is the set of positives with $|P_0|\le d$, and the testing outcome is $\vec{r}\in\{0,1\}^t$. We show that  Algorithm \ref{algorithm 1} will output $P_0$ given the input $\vec{r}$, that is, we show $P=P_0$.

Let $\F_0\subseteq\F$ be the set of column vectors of $M$ corresponding to $P_0$. Then, $|\F_0|\le d$ and $\vec{r}=\bigvee_{\vec{c}\in\F_0}\vec{c}$. According to Algorithm \ref{algorithm 1}, given the input $\vec{r}$, we first remove every column $\vec{c}\in\F$ that is not covered by $\vec{r}$. Then we obtain a subset $\F_S=\{\vec{c}_j: j\in S\}\subseteq\F$, where $S=\{j\in[n]:\vec{c}_j$ is covered by $\vec{r}\}$. It is obvious that $\F_0\subseteq\F_S$ and $\F_S\in\mathcal{U}(\F_0)$ where $\mathcal{U}(\F_0)$ is defined by (\ref{ssm-2}). Next we show how to determine $\F_0$ from $\F_S$.

For any $\vec{c}\in\F_S$, we claim that $\vec{c}\in\F_0$ if and only if there exists $i\in[t]$ such that $\vec{c}(i)=1$ and $\vec{c}'(i)=0$ for any $\vec{c}'\in\F_S\setminus\{\vec{c}\}$. To show the necessity, assume that $\vec{c}\in\F_0$ but there dose not exist $i\in[t]$ such that $\vec{c}(i)=1$ and $\vec{c}'(i)=0$ for any $\vec{c}'\in\F_S\setminus\{\vec{c}\}$. Then we have $\F_S\setminus\{\vec{c}\}\in\mathcal{U}(\F_0)$. Since $M$ is a $d$-SSM, we have $\vec{c}\not\in\bigcap_{\F'\in\mathcal{U}(\F_0)}\F'=\F_0$, a contradiction to the condition that $\vec{c}\in\F_0$. To show the sufficiency, assume that there exists $\vec{c}\in\F_S$ with the property that there exists $i_0\in[t]$ such that $\vec{c}(i_0)=1$ and $\vec{c}'(i_0)=0$ for any $\vec{c}'\in\F_S\setminus\{\vec{c}\}$, but $\vec{c}\not\in\F_0$, that is, the item corresponding to $\vec{c}$ is negative. Since $\F_0\subseteq\F_S$ and $\vec{r}=\bigvee_{\vec{c}\in\F_0}\vec{c}$, we must have $\vec{r}(i_0)=0$. Since $\vec{c}(i_0)=1$, it implies that $\vec{c}\not\in\F_S$ by the previous step of Algorithm \ref{algorithm 1}, a contradiction to the condition that $\vec{c}\in\F_S$.

Thus, the output $P$ of Algorithm \ref{algorithm 1} is exactly the set of positives $P_0$. The computational complexity of Algorithm \ref{algorithm 1} is $O(nt)$.
\end{IEEEproof}

\begin{algorithm}[H]\label{algorithm 1}
\caption{SSMIdAlg($\vec{r}$)}
Let $R_0=\{N_1,\ldots,N_{|R_0|}\}\subseteq[t]$ and $R_1=\{J_1,\ldots,J_{|R_1|}\}\subseteq[t]$ be two sets of indices which indicate $\vec{r}(i)=0$ and $\vec{r}(i)=1$ respectively. Clearly, $R_0\cup R_1=[t]$ and $|R_0|+|R_1|=t$.

$S=\{1,2,\ldots,n\}$;\\
$P=\emptyset$;\\
\For{$k=1$ to $|R_0|$}
 { $i=N_k$;\\
   \For{$j=1$ to $n$}
    { \If{$\vec{c}_j(i)=1$}
      { $S=S\setminus\{j\}$;
      }
    }
 }

\For{$k=1$ to $|R_1|$}
 { $i=J_k$;\\
   \For{$j=1$ to $n$}
    { \If{$j\in S$}
      { \If{$\vec{c}_j(i)=1$ and $\vec{c}_l(i)=0$ for any $l\in S\setminus\{j\}$}
        { $P=P\cup\{j\}$;
        }
      }
    }
 }
\eIf{$|P|\le d$}
 { {\bf output} $P$; }
 { {\bf output} ``The set of positives has size at least $d+1$.'' }
\end{algorithm}

\section{Bounds for $d$-SSM}
\label{sec:4}

In this section, we concentrate on the largest rate of a $d$-SSM. We first provide general bounds for $d$-SSM based on its connections with $d$-DM and $\bar{d}$-SM, and then derive an improved lower bound for $2$-SSM by the random coding method with expurgation, which is much better than the best existing lower bound of $2$-DM.

\subsection{General bounds for $d$-{\rm SSM}}
\label{sec:4-1}

Let $n(d,t)$, $f(d,t)$ and $s(\bar{d},t)$ denote the maximum possible number of columns of a $d$-SSM, a $d$-DM and a $\bar{d}$-SM with $t$ rows respectively. Denote their largest rates as
\begin{equation*}
  \begin{aligned}
  R(d)   & =\underset{t\to\infty}{\overline{\lim}}\frac{\log_2 n(d,t)}{t},\\
  R_D(d) & =\underset{t\to\infty}{\overline{\lim}}\frac{\log_2 f(d,t)}{t},\\
  R_S(\bar{d}) & =\underset{t\to\infty}{\overline{\lim}}\frac{\log_2 s(\bar{d},t)}{t}.
\end{aligned}
\end{equation*}

In the literature, the best known upper and lower bounds of $R_D(d)$ for $d\ge 3$ were proved in \cite{DR,DRR} respectively, and the general bounds of $R_S(\bar{d})$ were derived by Lemma \ref{relationship of DM and SM} and the known results on $R_D(d)$ \cite{DH}. For the case $d=2$, Erd\"{o}s, Frankl and F\"{u}redi \cite{EFF} investigated cover-free families and derived the best known results for $R_D(2)$ in which the lower bound was obtained by the random coding method and the upper bound was obtained by the techniques in extremal combinatorics. In \cite{CS}, Coppersmith and Shearer provided the best known lower bound for $R_S(\bar{2})$ by constructing a $\bar{2}$-union-free family from a deterministic \textit{cancellative} family and a random \textit{weakly union-free} family, and gave the best known upper bound for $R_S(\bar{2})$ by the techniques also in extremal combinatorics.

\begin{theorem}{\rm({\cite{CS,DH,DR,DRR,EFF}})}\label{known bounds of DM,SM}
Let $d\ge 2$ be an integer. If $d\to\infty$, then we have
\begin{equation*}
  \frac{1}{d^2\log_2e}(1+o(1)) \le R_D(d)\le R_S(\bar{d})\le R_D(d-1)\le \frac{2\log_2(d-1)}{(d-1)^2}(1+o(1))
\end{equation*}
where $e$ is the base of the natural logarithm. Moreover,
\begin{gather*}
  0.1814\le R_D(2) \le 0.3219, \\
  0.3135\le R_S(\bar{2})\le 0.4998.
\end{gather*}
\end{theorem}

We remark that the expressions of the general bounds of $R_D(d)$ for any $d\ge 3$ shown in \cite{DR,DRR} are complicated and therefore not stated in this paper. The interested reader may refer to the references therein. The asymptotic version on the bounds of $R_D(d)$ shown in Theorem \ref{known bounds of DM,SM} could also be found in \cite{D}.

By Theorem \ref{known bounds of DM,SM} and Lemma \ref{relationship2}, we immediately have the following results for SSM.

\begin{corollary}\label{general ssm}
Let $d\ge 2$ be an integer. Then we have
\begin{equation}\label{general bounds of SSM}
  \frac{1}{d^2\log_2e}(1+o(1)) \le R(d) \le \frac{2\log_2(d-1)}{(d-1)^2}(1+o(1))
\end{equation}
for $d\to\infty$, and
\begin{equation}\label{bounds of 2-SSM}
  0.1814 \le R(2) \le 0.4998.
\end{equation}
\end{corollary}

\subsection{An improved lower bound for $2$-{\rm SSM}}
\label{sec:4-2}

Inspired by the Kautz-Singleton construction for DM in \cite{KS} which is based on maximum distance separable codes and identity codes, in this part, we provide an improved lower bound of $R(2)$ by the random coding method together with a concatenated construction for SSM based on \textit{strongly separable codes}. For more applications of this method, the interested reader may refer to \cite{AS}.

\begin{theorem}\label{new lower bound of R(2)}
$R(2)\ge 0.2213$.
\end{theorem}

Before proving Theorem \ref{new lower bound of R(2)}, we do some preparations. Let $q\ge 2$ be an integer and $Q=\{0,1,\ldots,q-1\}$ be an alphabet. A set $\C=\{\vec{c}_1,\vec{c}_2,\ldots,\vec{c}_n\}\subseteq Q^t$ is called a $(t,n,q)$ code where each $\vec{c}_j$ is called a \textit{codeword}, $t$ is the \textit{length} of the code and $n$ is the \textit{code size}. For a code $\C\subseteq Q^t$, define the set of the $i$th coordinates of $\C$ as $\C(i)=\{\vec{c}(i)\in Q:\ \vec{c}=(\vec{c}(1),\vec{c}(2),\ldots,\vec{c}(t))\in\C\}$ for any $1\le i\le t$ and define the \textit{descendant code} of $\C$ as ${\rm desc}(\C)=\C(1)\times\C(2)\times\ldots\times\C(t).$
%

In \cite{JCM}, Jiang, Cheng and Miao introduced strongly separable codes in multimedia fingerprinting for the purpose of tracing back to all the traitors in an averaging collusion attack as efficiently as the well-known \textit{frameproof codes} but having a larger code size than frameproof codes.

\begin{definition}{\rm (\cite{JCM})}\label{ssc}
Let $\C$ be a $(t,n,q)$ code and $d\ge 2$ be an integer. $\C$ is called a \textit{strongly $\bar{d}$-separable code}, or briefly $\bar{d}$-SSC$(t,n,q)$, if for any $\C_0\subseteq\C$ with $1\le|\C_0|\le d$, we have
\begin{equation}\label{ssc-1}
  \bigcap\limits_{\C'\in \mathcal{S}(\C_0)}\C'=\C_0,
\end{equation}
where
\begin{equation}\label{ssc-2}
  \mathcal{S}(\C_0)=\{\C'\subseteq\C:\ {\rm desc}(\C')={\rm desc}(\C_0)\}.
\end{equation}
\end{definition}

\begin{remark}
An equivalent description of Definition \ref{ssc} is as follows. A $(t,n,q)$ code $\C$ is a $\bar{d}$-SSC if for any $\C_0,\C'\subseteq\C$ with $1\le |\C_0|\le d$,
$\C'\in\mathcal{S}(\C_0)$ implies that $\C_0\subseteq\C'$.
\end{remark}

\begin{remark}\label{ssc:2}
When $q=2$, an equivalent description of (\ref{ssc-2}) is that $\mathcal{S}(\C_0)=\{\C'\subseteq\C: \bigvee_{\vec{c}\in\C'}\vec{c}=\bigvee_{\vec{c}\in\C_0}\vec{c}$ and $\bigwedge_{\vec{c}\in\C'}\vec{c}=\bigwedge_{\vec{c}\in\C_0}\vec{c}\}.$
\end{remark}

We have the following observation on the relationship between SSM and SSC.

\begin{lemma}\label{lemma 1}
If there exists a $d$-SSM of size $t\times n$, then there exists a $\bar{d}$-SSC$(t,n,2)$.
\end{lemma}
\begin{IEEEproof}
Let $M$ be a $d$-SSM of size $t\times n$ and $\C$ be the set of column vectors of $M$. It is obvious that $\C$ is a $(t,n,2)$ code. By (\ref{ssm-2}) and Remark \ref{ssc:2}, for any $\C_0\subseteq\C$ with $1\le|\C_0|\le d$, we have $\mathcal{S}(\C_0)\subseteq \mathcal{U}(\C_0)$. Since $M$ is a $d$-SSM, we have $\C_0=\bigcap_{\C'\in \mathcal{U}(\C_0)}\C'\subseteq \bigcap_{\C'\in \mathcal{S}(\C_0)}\C'\subseteq\C_0$ which yields $\bigcap_{\C'\in \mathcal{S}(\C_0)}\C'=\C_0$. Thus $\C$ is a $\bar{d}$-SSC. 
\end{IEEEproof}

In \cite{JCM}, Jiang, Cheng and Miao also provided a concatenated construction for $\bar{d}$-SSC$(tq,n,2)$ based on $\bar{d}$-SSC$(t,n,q)$. We show that the $(tq,n,2)$ code they constructed is actually a $tq\times n$ $d$-SSM.

\begin{lemma}\label{lemma 2}
If there exists a $\bar{d}$-SSC$(t,n,q)$, then there exists a $d$-SSM of size $tq\times n$.
\end{lemma}
\begin{IEEEproof}
Let $\C=\{\vec{c}_1,\ldots,\vec{c}_n\}$ be a $\bar{d}$-SSC$(t,n,q)$ on $Q=\{0,1,\ldots,q-1\}$. For each $\vec{c}_j=(\vec{c}_j(1),\ldots,\vec{c}_j(t))$ $\in\C$, define $\vec{x}_j=(\vec{x}_j^1,\vec{x}_j^2,\ldots,\vec{x}_j^t)\in\{0,1\}^{tq}$ where
\begin{equation*}
  \vec{x}_j^i=(\vec{x}_j^i(0),\vec{x}_j^i(1),\ldots,\vec{x}_j^i(q-1))\in\{0,1\}^q
\end{equation*}
for any $1\le i\le t$ and
\begin{equation*}
  \vec{x}_j^i(k)=\left\{\begin{matrix}
           1, & \vec{c}_j(i)=k \\
           0, & $otherwise$
           \end{matrix}\right.
\end{equation*}
for any $0\le k\le q-1$. Let $\F=\{\vec{x}_1,\vec{x}_2,\ldots,\vec{x}_n\}$ be the set of column vectors of a matrix $M$. It is obvious that $M$ is a binary matrix of size $tq\times n$. We show that $M$ is a $d$-SSM.

For any $\F_0,\F'\subseteq\F$ with $|\F_0|\le d$, let $\C_0,\C'\subseteq\C$ denote the corresponding subsets of codewords to $\F_0,\F'$ respectively. Then $|\C_0|\le d$. If $\F'\in\mathcal{U}(\F_0)$, we must have $\C'\in \mathcal{S}(\C_0)$ according to the construction for $\F$. Since $\C$ is a $\bar{d}$-SSC, we have $\C_0\subseteq\C'$ yielding $\F_0\subseteq\F'$. Thus $M$ a $d$-SSM.
\end{IEEEproof}

To use Lemma \ref{lemma 2} to derive bounds of SSM, the results on $\bar{d}$-SSC with fixed $q$ and large $t$ is required. However, to the best of our knowledge, there is no known good result for this case in the literature. Therefore, in order to derive the lower bound of $R(2)$ in Theorem \ref{new lower bound of R(2)}, we shall first randomly construct a $\bar{2}$-SSC with fixed small $q$ and large $t$ and then exploit Lemma \ref{lemma 2} to obtain a $2$-SSM. To present the argument more precisely, we need the following concept of \textit{minimal frame}, which was also studied in \cite{JGC}.

\begin{definition}\label{frame}
Let $\C$ be a $(t,n,q)$ code. For any $\C_0,\C'\subseteq\C$, we call $\C'$ a \textit{frame} of $\C_0$ if $\C'\in\mathcal{S}(\C_0)$. Moreover, $\C'$ is called a \textit{minimal frame} of $\C_0$ if $\C'\in\mathcal{S}(\C_0)$ and $\C'\setminus\{\vec{c}\}\not\in\mathcal{S}(\C_0)$ for any $\vec{c}\in\C'$.
\end{definition}

\begin{lemma}\label{lemma 3}
Let $\C$ be a $(t,n,q)$ code. If $\C$ is not a $\bar{d}$-SSC, then there exist $\C_0\subseteq\C$ with $1\le |\C_0|\le d$ and a minimal frame $\C'$ of $\C_0$ such that $\C_0\not\subseteq\C'$.
\end{lemma}
\begin{IEEEproof}
If $\C$ is not a $\bar{d}$-SSC, then by Definition \ref{ssc}, there exist $\C_0\subseteq\C$ with $1\le|\C_0|\le d$ and a frame $\C'\subseteq\C$ of $\C_0$ such that $\C_0\not\subseteq\C'$. If $\C'$ is minimal, then it completes the proof. Otherwise, by Definition \ref{frame}, there must exist a codeword $\vec{c}\in\C'$ such that $\C'\setminus\{\vec{c}\}$ is still a frame of $\C_0$. Consider $\C'\setminus\{\vec{c}\}$ and repeat the process until it forms a minimal frame of $\C_0$. 
\end{IEEEproof}

The following result could be found in \cite{JGC} as well. For readers' convenience, we will give a self-contained proof of it.

\begin{lemma}\label{lemma 4}
Let $\C$ be a $(t,n,q)$ code. For any $\C_0\subseteq\C$ with $1\le|\C_0|\le d$, the minimal frame of $\C_0$ has size no more than $td-t+1$.
\end{lemma}
\begin{IEEEproof}
Suppose that $\C'$ is a minimal frame of $\C_0$. We count the number of codewords in $\C'$ by the order of coordinates. For the first coordinate, by Definition \ref{frame} and (\ref{ssc-2}), there exists $\C_1\subseteq\C'$ such that $\C_1(1)=\C_0(1)$ and $(\C_1\setminus\{\vec{c}\})(1)\neq\C_0(1)$ for any $\vec{c}\in\C_1$. It is obvious that $|\C_1|=|\C_0(1)|\le d$. Consider the second coordinate, then there exists $\C_2\subseteq\C'\setminus\C_1$ such that $(\C_2\cup\C_1)(2)=\C_0(2)$ and $(\C_2\cup\C_1\setminus\{\vec{c}\})(2)\neq\C_0(2)$ for any $\vec{c}\in\C_2$. Then we have $|\C_2|=|\C_0(2)|-|\C_1(2)|\le d-1$. Consider $\C'\setminus(\C_1\cup\C_2)$ and the $i$th coordinates in a similar way for $3\le i\le t$ until there will be no codewords left, which must occur since $\C'$ is a minimal frame of $\C_0$. Then we have $|\C'|\le d+(d-1)(t-1)=td-t+1$ as desired.
\end{IEEEproof}

Now we present the proof of Theorem \ref{new lower bound of R(2)}.

\begin{IEEEproof}[Proof of Theorem \ref{new lower bound of R(2)}]
Let $n>3$ and $\C=\{\vec{c}_1,\vec{c}_2,\ldots,\vec{c}_n\}$ be a collection of vectors of length $t$ where $\vec{c}_j=(\vec{c}_j(1), \vec{c}_j(2),\ldots,\vec{c}_j(t))$ and each $\vec{c}_j(i)$ is chosen uniformly and independently at random from a set $Q=\{0,1,\ldots,q-1\}$ with the probability that
\begin{equation*}
  {\rm Pr}(\vec{c}_j(i)=k)=1/q,\ \forall k\in Q
\end{equation*}
for any $1\le j\le n$ and $1\le i\le t$. The values of $t,n,q$ will be determined later.

For any vector $\vec{c}\in\C$, $\vec{c}$ is called {\it bad} if there exists $\C'=\{\vec{c}_0,\vec{c}_1,\ldots,\vec{c}_m\}\subseteq\C\setminus\{\vec{c}\}$ such that at least one of the following two cases occurs:
\begin{enumerate}[(1)]
\item there exists some $\vec{c}_i\in\C'$ such that ${\rm desc}(\C')={\rm desc}(\{\vec{c},\vec{c}_i\})$ with $1\le m\le t$;
\item there exists some $\vec{c}_i\in\C'$ such that ${\rm desc}(\C'\setminus\{\vec{c}_i\})={\rm desc}(\{\vec{c},\vec{c}_i\})$ with $2\le m\le t+1$.
\end{enumerate}

Let $\mathcal{T}$ be the collection of all bad vectors of $\C$ and $\widehat{\C}=\C\setminus\mathcal{T}$. We claim that all the vectors in $\widehat{\C}$ are distinct. If not, assume that there exist $\vec{a},\vec{b}\in\widehat{\C}$ such that $\vec{a}=\vec{b}$. Then, for any $\vec{x}\in\widehat{\C}\setminus\{\vec{a},\vec{b}\}$, we have ${\rm desc}(\{\vec{a},\vec{x}\})={\rm desc}(\{\vec{b},\vec{x}\})$ which implies that $\vec{a}$ is bad, a contradiction to the assumption. Hence, $\widehat{\C}$ is a $(t,n-|\mathcal{T}|,q)$ code by regarding each vector in $\widehat{\C}$ as a codeword.

We further show that $\widehat{\C}$ is a $\bar{2}$-SSC. If not, by Lemma \ref{lemma 3}, there exists $\C_0\subseteq\widehat{\C}$ with $1\le|\C_0|\le 2$ and a minimal frame $\C_1\subseteq\widehat{\C}$ of $\C_0$ such that $\C_0\not\subseteq\C_1$. We discuss the size of $\C_0$.
\begin{enumerate}[1)]
  \item If $|\C_0|=1$, then all the codewords in $\C_1$ are the same as that in $\C_0$. Since $\C_1$ is a minimal frame of $\C_0$, we have $|\C_1|=1$ and $\C_1=\C_0$, a contradiction to the fact that $\widehat{\C}$ is a $(t,n-|\mathcal{T}|,q)$ code.

  \item If $|\C_0|=2$, let $\C_0=\{\vec{c},\vec{c}_0\}$.
  \begin{enumerate}[2.1)]
  \item If $|\C_0\cap\C_1|=1$, without loss of generality, assume that $\C_0\cap\C_1=\{\vec{c}_0\}$ and $|\C_1|=m+1$. Then by Lemma \ref{lemma 4}, we have $m\le t$. If $m=0$, we have $\C_1=\{\vec{c}_0\}$ and ${\rm desc}(\{\vec{c},\vec{c}_0\})={\rm desc}(\{\vec{c}_0\})$, which implies that $\vec{c}=\vec{c}_0$, a contradiction. So, we have $1\le m\le t$. Then $\C_1$ satisfies case (1) implying that $\vec{c}$ is a bad codeword, a contradiction to the fact that $\widehat{\C}$ contains no bad codewords.

  \item If $|\C_0\cap\C_1|=0$, without loss of generality, assume that $\C_1=\{\vec{c}_1,\vec{c}_2,\ldots,\vec{c}_m\}\subseteq\widehat{\C}\setminus\C_0$. Then by Lemma \ref{lemma 4}, we have $m\le t+1$. If $m=1$, we have ${\rm desc}(\{\vec{c},\vec{c}_0\})={\rm desc}(\{\vec{c}_1\})$, which implies that $\vec{c}=\vec{c}_0=\vec{c}_1$, a contradiction. So, we have $2\le m\le t+1$. Then $\C_1$ satisfies case (2) implying that $\vec{c}$ is a bad codeword, also a contradiction.
  \end{enumerate}
\end{enumerate}
Thus, $\widehat{\C}$ is a $\bar{2}$-SSC$(t,n-|\mathcal{T}|,q)$ where $|\mathcal{T}|$ is a random variable due to the random construction of $\C$. Next we estimate the expected value of $|\mathcal{T}|$.

For any $\vec{c}\in\C$, let $\mathcal{S}_1(\vec{c})=\{\C'\subseteq\C\setminus\{\vec{c}\}:$ $|\C'|=m+1$ and $\C'$ satisfies case (1)$\}$ and $\mathcal{S}_2(\vec{c})=\{\C'\subseteq\C\setminus\{\vec{c}\}:$ $|\C'|=m+1$ and $\C'$ satisfies case (2)$\}$. Then, we have
\begin{align*}
{\rm Pr}(\vec{c}\ \rm is\ bad) &={\rm Pr}(|\mathcal{S}_1(\vec{c})|\ge 1\ {\rm or}\ |\mathcal{S}_2(\vec{c})|\ge 1)\\
 &\le {\rm Pr}(|\mathcal{S}_1(\vec{c})|\ge 1)+{\rm Pr}(|\mathcal{S}_2(\vec{c})|\ge 1)\\
 &\le {\rm E}(|\mathcal{S}_1(\vec{c})|)+{\rm E}(|\mathcal{S}_2(\vec{c})|) 
\end{align*}
where the last inequality is by Markov's inequality, and
\begin{equation}\label{expectation}
\begin{aligned}
 {\rm E}(|\mathcal{T}|) &=n\cdot {\rm Pr}(\vec{c}\ \rm is\ bad)\\
 &\le n\cdot ({\rm E}(|\mathcal{S}_1(\vec{c})|)+{\rm E}(|\mathcal{S}_2(\vec{c})|))
\end{aligned}
\end{equation}
where
\begin{align*}
{\rm E}(|\mathcal{S}_1(\vec{c})|) &=\sum_{m=1}^t\binom{n-1}{m+1}\binom{m+1}{1}{\rm Pr}({\rm desc}(\{\vec{c},\vec{c}_i\})={\rm desc}(\C'))\\
&=\sum_{m=1}^t\binom{n-1}{m+1}\binom{m+1}{1}\left((1/q)^{m+1}+(1-1/q)((2/q)^m-(1/q)^m)\right)^t\\
&\le \sum_{m=1}^t (m+1)n^{m+1}\left((2^m-1)q-(2^m-2)\right)^tq^{-(m+1)t}\\
&\le t\cdot \max_{1\le m\le t}\left\{(m+1)n^{m+1}\left((2^m-1)q-(2^m-2)\right)^tq^{-(m+1)t}\right\}
\end{align*}
and
\begin{align*}
{\rm E}(|\mathcal{S}_2(\vec{c})|) &=\sum_{m=2}^{t+1}\binom{n-1}{m+1}\binom{m+1}{1}{\rm Pr}\left({\rm desc}(\{\vec{c},\vec{c}_i\})={\rm desc}(\C'\setminus\{\vec{c}_i\})\right)\\
&=\sum_{m=2}^{t+1}\binom{n-1}{m+1}\binom{m+1}{1}\left((1/q)^{m+1}+(1-1/q)\left((2/q)^m-2(1/q)^m\right)\right)^t\\
&\le \sum_{m=2}^{t+1} (m+1)n^{m+1}\left((2^m-2)q-(2^m-3)\right)^tq^{-(m+1)t}\\
&\le t\cdot \max_{2\le m\le t+1}\left\{(m+1)n^{m+1}\left((2^m-2)q-(2^m-3)\right)^tq^{-(m+1)t}\right\}.
\end{align*}

For any $m\ge 1$ and $q\ge 2$, we have $(2^m-1)q-(2^m-2)>(2^m-2)q-(2^m-3)$. If
\begin{equation}\label{inequality 1}
 t\cdot \max\limits_{1\le m\le t+1}\left\{(m+1)n^{m+1}\left((2^m-1)q-(2^m-2)\right)^tq^{-(m+1)t}\right\}\le 1/3, 
\end{equation}
then ${{\rm E}}(|\mathcal{S}_1(\vec{c})|)\le1/3$ and ${{\rm E}}(|\mathcal{S}_2(\vec{c})|)\le 1/3$ for any $\vec{c}\in\C$. According to (\ref{expectation}), we have ${{\rm E}}(|\mathcal{T}|)\le 2n/3$, that is, the expected number of bad codewords in $\C$ is at most $2n/3$. By the random construction, there exists $\C$ such that it contains at most $2n/3$ bad codewords which implies that $\widehat{\C}$, obtained by deleting all the bad codewords from $\C$, is a $\bar{2}$-SSC with at least $n/3$ codewords. By Lemma \ref{lemma 2}, we can obtain a $2$-SSM with $tq$ rows and at least $n/3$ columns. Thus,
\begin{equation}\label{rate-1}
  R(2)\ge \underset{t\to\infty}{\overline{\lim}}\frac{\log_2(n/3)}{tq}=\underset{t\to\infty}{\overline{\lim}}\frac{\log_2 n}{tq}.
\end{equation}

Now we would like to maximize the lower bound of $R(2)$ in (\ref{rate-1}) under the restriction of (\ref{inequality 1}). It is obvious that (\ref{inequality 1}) is equivalent to that for any $1\le m\le t+1$,
\begin{equation}\label{inequality 2}
 t(m+1)n^{m+1}\left((2^m-1)q-(2^m-2)\right)^tq^{-(m+1)t}\le 1/3.
\end{equation}
By taking $\log_2$ on (\ref{inequality 2}) and do some simplifications, we have that for any $1\le m\le t+1$,
\begin{equation*}
 \frac{\log_2n}{tq}\le\frac{\log_2q}{q}-\frac{\log_2((2^m-1)q-(2^m-2))}{(m+1)q}-\frac{\log_2(3t(m+1))}{tq(m+1)}.
\end{equation*}
Take $t,n,q$ such that
\begin{equation*}
 \frac{\log_2n}{tq}=\frac{\log_2q}{q}-\max_{1\le m\le t+1}\left\{\frac{\log_2((2^m-1)q-(2^m-2))}{(m+1)q}-\frac{\log_2(3t(m+1))}{tq(m+1)}\right\}-\frac{\epsilon}{t}
\end{equation*}
where $\epsilon=o(t)$ is a real number. Then (\ref{inequality 1}) will be established and by (\ref{rate-1}), we have
\begin{align*}
 R(2) & \ge\underset{t\to\infty}{\overline{\lim}}\frac{\log_2 n}{tq} \\
      & =\frac{\log_2 q}{q}-\underset{t\to\infty}{\overline{\lim}}\max_{1\le m\le t+1}\left\{\frac{\log_2((2^m-1)q-(2^m-2))}{(m+1)q}-\frac{\log_2(3t(m+1))}{tq(m+1)}\right\} \\
      & \ge\frac{\log_2 q}{q}-\underset{t\to\infty}{\overline{\lim}}\max_{1\le m\le t+1}\left\{\frac{\log_2((2^m-1)q-(2^m-2))}{(m+1)q}\right\}
\end{align*}
for any $q\ge 2$. To make this lower bound as large as possible, we take $q=4$ and then
\begin{align*}
 R(2) & \ge \frac{1}{2}-\underset{t\to\infty}{\overline{\lim}}\max\limits_{1\le m\le t+1}\frac{\log_2((2^m-1)4-(2^m-2))}{4(m+1)} \\
      & =\frac{1}{2}-\frac{\log_2 22}{16} \\
      & \doteq 0.2213
\end{align*}
as desired.
\end{IEEEproof}

\section{Conclusion}
\label{sec:5}

In this paper, we introduced strongly separable matrices for nonadaptive CGT to identify a small set of positives from a large population. We showed that a $d$-SSM has weaker requirements than a $d$-DM, but provides an equally efficient identifying algorithm as a $d$-DM. A general bound on the largest rate of $d$-SSM was established. Besides, by the random coding method with expurgation, we derived an improved lower bound on the largest rate of $2$-SSM which is much higher than the best known result of $2$-DM. The results presented in this paper showed that a $2$-SSM could work better than a $2$-DM for nonadaptive CGT, which makes the research on SSM important. It is of interest to further improve the lower and upper bounds on the largest rate of $2$-SSM and explore the explicit constructions of optimal $2$-SSM. It is also interesting but more challenging to extend the argument for $2$-SSM to general $d$-SSM.


\ifCLASSOPTIONcaptionsoff
  \newpage
\fi

\end{document}